\documentclass[10pt]{article}

\usepackage{amscd,amsmath, amssymb, fancyhdr, mathbbol}
\usepackage[backref=page]{hyperref}

\renewcommand*{\backrefalt}[4]{%
	\ifcase #1 (Not cited.)%
	\or        (Cited on page~#2.)%
	\else      (Cited on pages~#2.)%
	\fi}

\hypersetup{
	colorlinks   = true,
	citecolor    = magenta,
	linkcolor    =blue,
	urlcolor     =magenta	
}

\numberwithin{equation}{section}

\newcommand{\version}{version 3.0,\ \ Sept 1, 2023}

\def\eqref#1{(\ref{#1})}
\newcommand{\goth}{\mathfrak}

\newcommand{\arrow}{{\:\longrightarrow\:}}
\newcommand{\Z}{{\Bbb Z}}

\newcommand{\R}{{\Bbb R}}
\newcommand{\Q}{{\Bbb Q}}

\def\1{\sqrt{-1}\:}
\newcommand{\restrict}[1]{{\left|_{{\phantom{|}\!\!}_{#1}}\right.}}
\newcommand{\cntrct}                
{\hspace{2pt}\raisebox{1pt}{\text{$\lrcorner$}}\hspace{2pt}}

\renewcommand{\phi}{\varphi}
\renewcommand{\epsilon}{\varepsilon}
\renewcommand{\geq}{\geqslant}
\renewcommand{\leq}{\leqslant}
\renewcommand{\max}{{\rm max}}

\newcommand{\Teich}{\operatorname{Teich}}

\newcommand{\End}{\operatorname{End}}

\newcommand{\Pos}{\operatorname{Pos}}
\newcommand{\Pic}{\operatorname{Pic}}

\newcommand{\Aut}{\operatorname{Aut}}

\newcommand{\Diff}{\operatorname{Diff}}

\newcommand{\coker}{\operatorname{coker}}

\newcommand{\Kah}{\operatorname{Kah}}

\newcommand{\Comp}{\operatorname{Comp}}

\newcounter{Mycounter}[section]
\newcounter{lemma}[section]
\newcounter{claim}[section]
\newcounter{sublemma}[section]
\newcounter{corollary}[section]
\renewcommand{\thecorollary}{{Corollary \thesection.\arabic{corollary}}}
\newcommand{\corollary}{%
    \setcounter{corollary}{\value{Mycounter}}
    \refstepcounter{corollary}
    \stepcounter{Mycounter}
    {\noindent \bf \thecorollary:\ }}

\newcounter{theorem}[section]
\renewcommand{\thetheorem}{{Theorem \thesection.\arabic{theorem}}}
\newcommand{\theorem}{%
    \setcounter{theorem}{\value{Mycounter}}
    \refstepcounter{theorem}
    \stepcounter{Mycounter}
    {\noindent \bf \thetheorem:\ }}

\newcounter{conjecture}[section]
\newcounter{proposition}[section]
\renewcommand{\theproposition}
      {{Proposition \thesection.\arabic{proposition}}}
\newcommand{\proposition}{%
    \setcounter{proposition}{\value{Mycounter}}
    \refstepcounter{proposition}
    \stepcounter{Mycounter}
    {\noindent \bf \theproposition:\ }}

\newcounter{definition}[section]
\renewcommand{\thedefinition}
      {{Definition~\thesection.\arabic{definition}}}
\newcommand{\definition}{%
    \setcounter{definition}{\value{Mycounter}}
    \refstepcounter{definition}
    \stepcounter{Mycounter}
    {\noindent \bf \thedefinition:\ }}

\newcounter{example}[section]
\newcounter{remark}[section]
\renewcommand{\theremark}{{Remark \thesection.\arabic{remark}}}
\newcommand{\remark}{%
    \setcounter{remark}{\value{Mycounter}}
    \refstepcounter{remark}
    \stepcounter{Mycounter}
    {\noindent \bf \theremark:\ }}

\newcounter{problem}[section]
\newcounter{question}[section]
\newcommand{\proof}{\noindent{\bf Proof:\ }}

\newcommand{\pstep}{{\bf Proof. Step 1:\ }}

\makeatletter

\@addtoreset{equation}{section} \@addtoreset{footnote}{section}
\makeatother

\def\blacksquare{\hbox{\vrule width 5pt height 5pt depth 0pt}}
\def\endproof{\blacksquare}

\begin{document}
\begin{center}
{\LARGE\bf
Roundness of the ample cone and existence of double 
Lagrangian fibrations on hyperk\"ahler manifolds\\[4mm]
}

Ljudmila Kamenova\footnote{Partially supported 
by a grant from the Simons Foundation/SFARI (522730, LK)}, 
Misha Verbitsky\footnote{Partially supported 
				by the HSE University Basic Research Program,
FAPERJ 	SEI-260003/000410/2023 and CNPq - Process 310952/2021-2.
}

\end{center}

{\small \hspace{0.10\linewidth}
\begin{minipage}[t]{0.85\linewidth}
{\bf Abstract.} 
Let $M$ be a hyperk\"ahler manifold of maximal holonomy
(that is, an IHS manifold), and let $K$ be its K\"ahler cone,
which is an open, convex subset in the space 
$H^{1,1}(M, \R)$ of real (1,1)-forms. 
This space is equipped with a
canonical bilinear symmetric form of signature $(1,n)$
obtained as a restriction of the 
Bogomolov-Beauville-Fujiki form. The set of
vectors of positive square in the space of signature $(1,n)$
is a disconnected union of two convex cones. The ``positive
cone" is the component which contains the K\"ahler cone. We
say that the K\"ahler cone is ``round" if it is equal to the
positive cone. The manifolds with round K\"ahler cones have
a unique bimeromorphic model and correspond to Hausdorff
points in the corresponding Teichm\"uller space. 
We prove that any maximal holonomy hyperk\"ahler manifold
with $b_2 > 5$ has a deformation with round K\"ahler cone
and the Picard lattice of signature (1,1), admitting two
non-collinear integer isotropic classes. This is used to 
show that all known examples of hyperk\"ahler manifolds admit a
deformation with two transversal Lagrangian fibrations,
and their Kobayashi metric vanishes, unless the
Picard rank is maximal.
\end{minipage}
}


\section{Introduction} 

This paper gives a simple solution
for a construction problem of hyperk\"ahler
geometry. We construct a hyperk\"ahler manifold with
rank 2 Picard lattice of signature (1,1) containing 
isotropic vectors and maximal
possible K\"ahler cone. 

For our purposes, ``a hyperk\"ahler manifold''
is a compact, holomorphic symplectic compact
manifold $M$ of K\"ahler type, which satisfies
``the maximal holonomy condition'', that is,
$\pi_1(M)=0$, $\dim H^{2,0}(M)=1$.
This condition is also known as IHS 
(``irreducible holomorphic symplectic'').

The shape of the K\"ahler cone of a hyperk\"ahler
manifold is more or less understood by now 
(see \cite{_AV:hyperbolic_cusps_}). However,
finding examples of manifolds with prescribed
shape of their K\"ahler cone is a
complicated task. Such constructions are 
either very explicit or based on convoluted
arguments from number theory. The automorphism
group of a hyperk\"ahler manifold and its set of
Lagrangian fibrations can be
described explicitly in terms of its periods and the shape of the
K\"ahler cone 
(\cite{_AV:construction_autom_}). Therefore, finding manifolds with prescribed
K\"ahler cones has many practical applications.

Recall that the second cohomology of a maximal holonomy
hyperk\"ahler manifold is equipped with a 
bilinear symmetric form of signature $(3, b_2-3)$,
which is essentially of topological origin
(\cite{_Beauville_,_Bogomolov:defo_,_Fujiki:HK_}).
This form, denoted by $q$ further on, 
is called {\bf the Bogomolov-Beauville-Fujiki
  form}; it is positive on the K\"ahler cone, and
has signature $(1,b_2-3)$ on $H^{1,1}(M,\R)$.

This implies that the set 
of ``positive vectors'' (the vectors with positive
square) in $H^{1,1}(M,\R)$ has two connected
components, both of them convex cones. However,
only one of these two components may contain K\"ahler
forms. We call this component ``the positive cone''
of a hyperk\"ahler manifold.

Every  hyperk\"ahler manifold is equipped
with a collection ${\goth S}$ of primitive integer 
cohomology classes in $ H^2(M, \Z)$ with negative squares,
called {\bf the MBM classes} (\cite{_AV:MBM_}).
This set is invariant on deformations of $M$
and under the action of the
monodromy group $\Gamma$ generated by the
monodromy of the Gauss-Manin connection
for all deformations of $M$. The
group $\Gamma$ (originally defined by E. Markman,
\cite{_Markman:mono_,_Markman:survey_}, who called it {\bf the monodromy group})
is mapped to the orthogonal lattive
 $O(H^2(M, \Z))$ with finite kernel,
and its image is a finite index sublattice
in $O(H^2(M, \Z))$ (\cite{_V:Torelli_}).

In \cite{_AV:MK_}, it was shown that the
monodromy group $\Gamma$
acts on the set of MBM classes
with a finite number of orbits,
which were computed for some
deformational classes of hyperk\"ahler manifolds
in \cite{_Bayer_Macri_},
\cite{_HT:intersection_}, \cite{_HT:Extremal_},
\cite{_AV:MBM-examples_}.

An {\bf MBM bound} for a hyperk\"ahler manifold is the
number \[ C:= \max \{ -q(x,x)\ \ |\ \ x\in {\goth S}\},\]
where ${\goth S}$ denotes the set of MBM classes.
Since $\Gamma\subset O(H^2(M, \Z))$ acts
on ${\goth S}$ with finitely many orbits
(by \cite{_AV:MK_}), this number is finite.

As shown in \cite{_AV:hyperbolic_cusps_}
the positive cone $\Pos(M,I)$ of a hyperk\"ahler manifold
is cut into pieces by hyperplanes orthogonal to
the MBM classes which lie in $H^{1,1}(M,I)$, and
each of the connected components of this complement
can be realized as a K\"ahler cone of a certain
hyperk\"ahler birational model of $(M,I)$.
In other words, the K\"ahler cone is
a connected component of the set
\begin{equation}\label{_Kahler_cone_split_Equation_}
\Pos^\circ(M,I):= \Pos(M,I) \left\backslash \bigcup_{\alpha \in
  {\goth S}\cap H^{1,1}(M,I)} \alpha^\bot\right.,
\end{equation}
where ${\goth S}$ is the set of all MBM classes in $H^2(M, \Z)$,
and all connected components are realized as K\"ahler cones
for birational models of $(M,I)$

The authomorphism group of a hyperk\"ahler
manifold $(M,I)$ is expressed in terms of its K\"ahler 
cone and the monodromy group as follows.
Let  $\Gamma_I\subset \Gamma$ be the subgroup of
the monodromy group preserving the Hodge
decomposition on $H^2(M)$. Then $\Aut(M,I)$
is a subgroup of all elements in $\Gamma_I$ which 
preserve the K\"ahler cone
(the result is essentially due to
E. Markman,  \cite{_Markman:survey_}),

We say that a manifold $M$ {\bf has round K\"ahler cone}
if $\Kah(M)=\Pos(M)$, or, equivalently, 
when the set of MBM classes in $H^{1,1}(M,I)$ is empty.

\hfill

Our main results are the following two theorems.

\hfill

\hfill

\noindent
{\bf \ref{_vanishing_main_Theorem_}:\;}
Let $M$ be a compact maximal holonomy hyperk\"ahler
manifold  with $b_2(M)\geq 4$, satisfying the SYZ 
conjecture. Assume that 
$H^2(M,\Q)$ has non-zero isotropic vectors.
Then $M$ admits a deformation 
with two distinct Lagrangian fibrations.
If, in addition, $M$ 
satisfies one of the two assumptions
\begin{description}
\item[(a)] 
$(H^{2,0}(M) \oplus H^{0,2}(M)) \cap H^2(M, \Q)=0$,  
\item[(b)]  $b_2(M)\geq 6$, and $M$ has 
Picard lattice of non-maximal rank, 
\end{description}
then the Kobayashi pseudometric on $M$ vanishes.

\hfill

\noindent
{\bf \ref{_vanishing_KObayashi_Corollary_}:\;}
All known compact hyperk\"ahler examples have vanishing Kobayashi 
pseudometric.


\section{Integral quadratic lattices of rank 2}



%

Recall that 
a sublattice $R \subset \Z^n$ is called {\bf primitive}
if $(R\otimes_\Z \Q) \cap \Z^n=R$. The {\bf saturation}
of a sublattice $R\subset \Z^n$ is the sublattice
$R_1:= (R\otimes_\Z \Q) \cap \Z^n$. Clearly,
$R_1\supset R$ is a lattice of the same rank.


\hfill

We say that an integral 
quadratic lattice $(\Lambda, q)$ {\bf represents $n\in \Z$}
if $q(x,x)=n$ for some $x\in \Lambda$.
In this paper, we prove the following technical
result, used to obtain deformations of a hyperk\"ahler
manifold without MBM classes.

\hfill

\proposition \label{_1,1_sublattice_with_MBM_bound_Proposition_}
Let $\Lambda$ be a quadratic integral lattice of signature 
$(3, k)$, $k \geq 1$, representing 0. 
Then for any $N>0$ the 
lattice $\Lambda$ contains a primitive
integral sublattice of signature $(1,1)$
representing $0$ and not representing integers
in the interval $[-N, -1]$.

\hfill

\pstep
To find a  a primitive
integral sublattice of signature $(1,1)$ {\em not}
representing $0$ and not representing integers
in the interval $[-N, -1]$, we could apply
\cite[Theorem 3.6]{_AV:construction_autom_}.
Now we need to prove the same result for
lattices representing 0.

Choose an intermediate lattice
$B\subset \Lambda$, not necessarily primitive,
and let $B^s$ be its saturation. Let $d$ be
the index of $B$ in $B^s$.
For any sublattice $L\subset B$ primitive
in $B$, denote by $L^s$ the saturation of $L$ in $\Lambda$.
Consider the diagram
\[
\begin{CD}
0 @>>> L @>>>  L^s @>>>  L^s/L @>>> 0\\
& & @V\nu VV @V{\nu_s}VV @VV {\nu_c} V \\
0 @>>> B @>>>  B^s @>>>  B^s/B @>>> 0
\end{CD}
\]
The first two vertical arrows $\nu$ and $\nu_s$  are injective by 
construction.
The Snake Lemma implies the long exact sequence 
\[
0 \arrow \ker \nu_c \stackrel \delta \arrow \coker \nu\arrow 
\coker \nu_s\arrow \dots
\]
Since $L$ is primitive in $B$, then $\coker \nu_c$ is 
torsion-free. However, $L^s/L$ is a torsion group 
by construction. Therefore, the map $\delta$ 
in this exact sequence vanishes, which implies
the vanishing of $\ker \nu_c$. We obtain that
$L^s/L$ is embedded into $B^s/B$, and 
$|L^s/L|\leq R:=|B^s/B|$. 

Assume that the quadratic form 
$q\restrict L$ does not take any values on 
the interval $[-N, -1]$. Since all
elements of $L^s$ are proportional to
elements in $L$, with the coefficient
of proportionality bounded by $|L^s/L|$,
this implies that $q\restrict{L^s}$
does not take any of the values in  $[-\frac N{|L^s/L|}, -1]$.
Therefore, $|L^s/L|\leq R$ implies that
$q\restrict{L^s}$
does not take values in $[-\frac N {R^2}, -1]$.

Any non-degenerate lattice over $\Q$ has
an orthogonal basis. If it represents 0 and 
its signature is $(3, k), k >0$, it also 
has a basis $(x, y, z_1, ..., z_{k+1})$
such that $x$ and $y$ are isotropic,
$q(x, y) =1$, and $z_i$ are pairwise 
orthogonal and orthogonal to $x$ and $y$.
The corresponding Gram matrix has the form
\[
\begin{pmatrix}
0 & 1 & 0 & \cdots & 0\\
1 & 0 & 0 & \cdots & 0\\
0 & 0 & \alpha_1 & \cdots & 0\\
\vdots & \vdots &\vdots  & \ddots & \vdots\\
0 & 0 & 0& \cdots & \alpha_{k+1}
\end{pmatrix}
\] 
we will call a quadratic form with such a
matrix {\bf hyperbolic-orthogonal}. 

Let $B\subset \Lambda$ be 
generated by integral vectors proportional
to such a basis in $\Lambda\otimes \Q$. To prove
\ref{_1,1_sublattice_with_MBM_bound_Proposition_},
it remains to show that any hyperbolic-orthogonal
quadratic integral lattice $B$ of signature 
$(3, k)$, $k \geq 1$ contains a primitive signature 
(1,1) sublattice representing 0 and
not representing any integers in the interval $[-R^2N, -1]$.

Let $v= \sum_{i=1}^{k+1} \lambda_i z_i$ be an integral
vector in $B$; it is primitive if and only if
the greatest common divisor of $\lambda_i$ is 1.
A vector $x + t (y+v)$ is isotropic if and only if
\[ 
0 =2tq(x,v+y)+ t^2 q(v+y,v+y)=2t + t^2 q(v,v)
\]
or, equivalently, if $t= -\frac{2}{q(v,v)}$.
This gives a isotropic integral vector
$q(v,v)x - 2 (v+y)$, which is primitive
whenever $q(v,v)$ is odd. When $q(v,v)$ is even,
the vector $\frac{q(v,v)}2 x -v -y$ is primitive
and isotropic. Denote this isotropic vector
by $x'$. Then the rank 2 lattice $\langle y, x'\rangle$
is primitive in the lattice $B$ defined above.
Indeed, $\alpha y + \beta x'= \alpha y + \beta (q(v,v)x - 2 v)$
is primitive when $\alpha$ and $\beta$ are coprime,
because $x$ and $y$ are parts of the base of $B$.

It remains to show that $\langle y, x'\rangle$,
for an appropriate choice of $v$, does not
represent integers in the interval $[-R^2 N, -1]$.

Clearly, $q(y, x')$ is either $q(v,v)$
or $\frac{q(v,v)}2$, depending on parity of $q(v,v)$.
A general vector in the integral lattice $L:= \langle y, x'\rangle$
has the form $\beta (q(v,v)x - 2 v) + \alpha y$ or
$\beta (\frac{q(v,v)} 2 x -  v) + \alpha y$
depending on parity of $q(v,v)$, where $\alpha$ and
$\beta$ are arbitrary integers. The square of this
vector is $4\beta^2 q(v,v) + 2q(v,v)\alpha \beta$
in the first case, and  $\beta^2 q(v,v)^2 + \frac {q(v,v)} 2 \alpha \beta$
in the second case. In both cases, $L$ does not represent
any integer $t \in (-\frac{q(v,v)} 2 , \frac{q(v,v)} 2) \backslash 0$.
Choosing a vector $v$ with sufficiently big square, we
can make sure that the bound $\frac{q(v,v)} 2 > R^2N$
is reached, and the lattice $L^s$
does not represent any integer in 
the interval $[-N, -1]$.
\endproof


\section{Kobayashi metric on hyperk\"ahler manifolds}

We apply the results of the current paper
to the vanishing of the Kobayashi pseudometric on hyperk\"ahler
manifolds. In this section we summarize some of our results 
in \cite{klv} joint with S. Lu. The aim of the current paper is to 
imrove some of the bounds imposed on the Betti numbers, and also 
to show vanishing of the Kobayashi pseudometric for all of the known 
compact hyperk\"ahler examples. This generalizes our result on Kobayashi 
non-hyperbolicity of all known examples, \cite{_Kamenova_V:fibrations_}. 

\hfill

\definition 
An {\bf ergodic complex structure}
is a complex structure $I$ on $M$ such that 
for any complex structure $I'$ in the same
deformation class there exists a sequence
of diffeomorphisms $\nu_i\in \Diff(M)$
such that $\lim_i \nu_i(I)=I'$,
where the limit of $\nu_i(I)\in \End(TM)$
is taken with respect to the $C^\infty$-topology on the
space of tensors. We denote the space of
all integrable complex structures with
this topology by $\Comp$.

\hfill

\theorem\label{_ergodic_big_orbit_Theorem_}
Any complex structure of hyperk\"ahler type 
on a hyperk\"ahler manifold with $b_2 \geq 5$
with $(H^{2,0}(M)\oplus H^{0,2}(M))\cap H^2(M, \Q)=0$
is ergodic.

\proof
\cite{_Verbitsky:ergodic_,_Verbitsky:ergodic_erratum_}).
\endproof

\hfill

We will need another diffeomorphism orbit, which is
smaller than the maximal one, but has 
many of the same properties.

\hfill

\theorem\label{_smaller_orbit_ergodic_Theorem_}
Let $(M,I)$ be a hyperk\"ahler manifold
such that $(H^{2,0}(M)\oplus H^{0,2}(M))\cap H^2(M, \Q)$
is a rank 1 space generated by a class  
$\alpha\in  H^2(M, \Q)$, and $\Teich_\alpha$ the
Teichm\"uller space of all complex structure
with $\alpha\in (H^{2,0}(M)\oplus H^{0,2}(M))\cap H^2(M, \Q)$
and deformationally equivalent to $I$. Then 
$\Diff(M) \cdot I$ is dense in $\Teich_\alpha$.

\proof
\cite[Theorem 2.5, Theorem 3.1]{_Verbitsky:ergodic_erratum_}.
\endproof

\hfill

\theorem\label{_ergodic_vanishes_Theorem_} 
\cite{klv}
Let $(M,I)$ be a complex manifold 
with vanishing Kobayashi pseudometric. Then the Kobayashi 
pseudometric vanishes for all ergodic complex structures 
in the same deformation class. Moreover, for each
complex structure $I_1$ such that the closure of 
$\Diff(M) I_1$ in $\Comp$ contains $I$,
the pseudometric on $(M,I_1)$ also vanishes.

\hfill

\proof The proof follows easily from semicontinuity
of the diameter of the Kobayashi pseudometric,
considered as a function on $\Comp$ (\cite{klv}). \endproof

\hfill

\theorem\label{_two_Lagra_vanishing_Theorem_} \cite{klv} 
Let $M$ be a hyperk\"ahler manifold admitting two 
Lagrangian fibrations associated with two non-proportional 
parabolic classes. Then the Kobayashi pseudometric on $M$ vanishes. 

\hfill

To prove that a given hyperk\"ahler manifold admits
a deformation with two distinct Lagrangian fibrations, in \cite{klv} 
we used an argument based on \cite{_AV:construction_autom_}.

\hfill

\theorem\label{_manifold_with_parabo_AV_Theorem_} \cite{klv}
Let $M$ be a maximal holonomy hyperk\"ahler manifold
with $b_2(M) > 13$. Then $M$ admits a projective deformation
with Picard lattice of signature $(1,2)$, with 
round K\"ahler cone (that is, with the K\"ahler cone
equal to the positive cone), and its automorphism group
has finite index in the arithmetic group $SO(\Pic(M))$
of orthogonal automorphisms of its Picard lattice.

\hfill

{\bf Proof:} From \cite[Theorem 3.11]{_AV:construction_autom_}
it follows that there exists a projective deformation with
Picard rank 3, isotropic classes in $H^{1,1}(M)\cap H^2(M, \Q)$ 
and without MBM classes of type (1,1). 
From \cite[Theorem 2.10]{_AV:construction_autom_}
it follows that for such a manifold the  K\"ahler cone
is equal to the positive cone, and from 
\cite[Theorem 2.6, Theorem 2.7, Corollary 2.12]{_AV:construction_autom_}
it follows that its automorphism group
has finite index in the arithmetic group $SO(\Pic(M))$.
\endproof

\hfill

\theorem\label{_exi_defo_2_fibra_Theorem_} \cite{klv} 
Let $M$ be a projective, maximal holonomy hyperk\"ahler manifold
with two non-collinear isotropic rational classes
in $H^{1,1}(M)$ and with round K\"ahler cone. 
Assume that $M$ satisfies the
SYZ conjecture, that is, any nef bundle on
$M$ is semiample. Then $M$ admits at least two 
transversal holomorphic Lagrangian fibrations.
In particular, the Kobayashi pseudometric on $M$ vanishes. 

\hfill

{\bf Proof:} Since the K\"ahler cone of $M$
is round,  there exist rational vectors on the boundary of the
K\"ahler cone of $M$.
These points correspond to rational points in the
real quadric $\{l \in {\Bbb P} H^2(M, \Q)\ \ |\ \ q(l,l)=0\}$.

Each of such points corresponds to a Lagrangian fibration,
because we assume that the SYZ conjecture holds. 
\endproof

\hfill

\theorem \cite{klv}
Let $(M,I)$ be a maximal holonomy, compact hyperk\"ahler manifold
with non-maximal Picard rank.
Suppose that it has a deformation which has 
two transversal Lagrangian fibrations.
Then the Kobayashi pseudometric on $(M,I)$ vanishes.

\hfill

{\bf Proof:} 
The vanishing of the Kobayashi pseudometric
then immediately follows from 
\ref{_exi_defo_2_fibra_Theorem_} and
\ref{_ergodic_vanishes_Theorem_}.
Indeed, there exists a complex structure 
$I'$ with vanishing Kobayashi pseudometric,
and a sequence of diffeomorphisms
such that $\lim_i \nu_i(I)=I'$.
Then the Kobayashi pseudometric of $(M,I)$
vanishes by semicontinuity properties of the diameter of the 
Kobayashi pseudometric.
 \endproof

\hfill

\theorem 
Let $M$ be a compact, maximal holonomy hyperk\"ahler manifold with 
$b_2(M)\geq 6$. Suppose that all its deformations satisfy 
the SYZ conjecture.  Then the Kobayashi pseudometric on $M$ vanishes.

\hfill

{\bf Proof:} See
Consider a rank 2 primitive sublattice 
$L\subset H^2(M, \Z)$ not representing
any integers in $[-N, -1]$, for
$N$ sufficiently big. By the global Torelli theorem
(\cite{_V:Torelli_}), for any K3-type Hodge structure
on $H^2(M,\R)$ there exists a complex structure $I$ of hyperk\"ahler
type on $M$ with the same Hodge decomposition.
Let $I_L$ be a complex structure such that $H^{1,1}_{I_L}(M,\Z)= L$
and $\langle H^{2,0}_{I_L}(M) +  H^{0,2}_{I_L}(M)\rangle \cap H^{2}(M,\Z)=0$.
Then \ref{_ergodic_big_orbit_Theorem_} implies that $I_L$ is ergodic.
Using \ref{_exi_defo_2_fibra_Theorem_}, 
\ref{_ergodic_vanishes_Theorem_} and
\ref{_two_Lagra_vanishing_Theorem_}, we obtain
that the Kobayashi pseudometric vanishes on $(M,I)$
for any complex structure $I$ of hyperk\"ahler type in the same
deformation class as $I_L$. \endproof

\hfill

\remark
All known examples of hyperk\"ahler manifolds 
have $b_2(M)\geq 7$ and satisfy the SYZ conjecture.
By the results above, the Kobayashi
pseudometric of all known manifolds vanishes,
unless their Picard rank is maximal.

\hfill

\remark \label{SYZ_ref}
The SYZ conjecture is true for all known hyperk\"ahler examples.
Using the Fourier-Mukai transform and the deformations to
moduli spaces, it was proven, in projective case,
for the deformations of Hilbert schemes of points on K3 surfaces 
(Bayer-Macr\`i \cite[Theorem 1.5]{_Bayer_Macri_}; Markman 
\cite[Theorems 1.3 and 6.3]{_Markman:SYZ_}), 
for the deformations of the generalized Kummer varieties (Yoshioka 
\cite[Proposition 3.38]{_Yoshioka_}), for the O'Grady's sixfolds 
(Mongardi-Rapagnetta \cite[Corollary 1.3 and 7.3]{_Mongardi_Rapagnetta_}), 
and the for O'Grady's tenfolds (Mongardi-Onorati, 
\cite[Theorem 2.2]{_Mongardi_Onorati_}). 
In fact, in all these references except \cite{_Yoshioka_},
the result is stated for all K\"ahler deformations.
However, the passage from SYZ conjecture for 
all projective deformations to SYZ for all K\"ahler
is known in general, as shown in 
\cite[Theorem 1.2]{_Matsushita:Isotropic_} 
(see also 
\cite[Theorem 3.4]{_Kamenova_V:fibrations_}).
\endproof


\section{Main results}

The main result of this note is the following. 

\hfill

\theorem\label{_vanishing_main_Theorem_}
Let $M$ be a compact maximal holonomy hyperk\"ahler
manifold  with $b_2(M)\geq 4$, satisfying the SYZ 
conjecture. Assume that 
$H^2(M,\Q)$ has non-zero isotropic vectors.
Then $M$ admits a deformation 
with two distinct Lagrangian fibrations.
If, in addition, $M$ 
satisfies one of the two assumptions
\begin{description}
\item[(a)] 
$(H^{2,0}(M) \oplus H^{0,2}(M)) \cap H^2(M, \Q)=0$,  
\item[(b)]  $b_2(M)\geq 6$, and $M$ has 
Picard lattice of non-maximal rank, 
\end{description}
then the Kobayashi pseudometric on $M$ vanishes.

\hfill

{\bf Proof:}
Consider a primitive lattice $\Lambda \subset H^2(M, \Z)$
of signature $(p, q)$, $p \leq 1, q \leq b_2 -3$. 
Using the global Torelli theorem \cite{_V:Torelli_},
we can find a deformation $(M,I_1)$ of $M$ with Picard
lattice $\Lambda$.

Since the rank of the indefinite lattice $H^2(M,\Z)$ is at least $5$, 
by Meyer's theorem \cite{_Meyers_} there exists an isotropic
 vector $x \in H^2(M, \Z)$. Applying
 \ref{_1,1_sublattice_with_MBM_bound_Proposition_}
we can find a primitive integral sublattice $L^s\subset H^2(M,\Z)$  
of signature (1,1) representing 0 and not representing
all integers in $[-C, -1]$, where $C$ is the MBM bound.

Choose the complex structure
$I$ such that $L^s$ is the Picard lattice
$H^{1,1}_I(M,\Z)$ of $(M,I)$. Then
$(M,I)$ has round K\"ahler cone, 
hence both integer isotropic generators
of $U_N$ are nef. 

These classes correspond
to Lagrangian fibrations since $(M,I)$ 
satisfies SYZ, hence the Kobayashi metric
of $(M,I)$ vanishes. This takes care of the
first statement of \ref{_vanishing_main_Theorem_}.

Applying \ref{_ergodic_vanishes_Theorem_},
we obtain that the Kobayashi metric vanishes
for all ergodic complex structures,
that is, for all complex structures $I_1$
such that $(H^{2,0}(M) \oplus H^{0,2}(M)) \cap H^2(M,
\Q)=0$. This proves the case (a) of
 \ref{_vanishing_main_Theorem_}.

It remains to prove  \ref{_vanishing_main_Theorem_} (b).
Suppose that 
\[ \bigg(H^{2,0}(M,I) \oplus H^{0,2}(M,I)\bigg) \cap H^2(M, \Q)
\]
has rank one and is generated by $\alpha$.
Since $\alpha^\bot$ has rank $\geq 5$,
it contains isotropic vectors. Applying
\ref{_1,1_sublattice_with_MBM_bound_Proposition_} again, we find 
a deformation $(M,I')$ of $M$ which satisfies
$H^{1,1}_{I'}(M,\Z) = L^s$ and
$(H^{2,0}(M,I') \oplus H^{0,2}(M,I')) \cap H^2(M, \Q)=\langle\alpha\rangle$.
For an appropriate choice of diffeomorphisms
$\pi_i \in \Diff(M)$, the sequence $\nu_i(I')$ 
converges to $I$ (by \ref{_smaller_orbit_ergodic_Theorem_}),
hence the Kobayashi metric on $(M,I)$ also vanishes.
\endproof

\hfill

\corollary \label{_vanishing_KObayashi_Corollary_}
All known compact hyperk\"ahler examples have vanishing Koba\-yashi 
pseudometric. 

\hfill

{\bf Proof:} See \ref{SYZ_ref} and \ref{_vanishing_main_Theorem_}. 
\endproof

\hfill

{\bf Acknowledgments.} 
We are grateful to Christian Lehn for pointing out a gap in an 
argument from our paper \cite{klv}. We also thank Giovanni Mongardi for 
his suggestions and references. The referee was very helpful with their 
comments, improving the paper and its exposition, and finding
an important error in the exposition. We are grateful to Ivan
Frolov, who made it possible to fix this error.

\noindent {\sc Ljudmila Kamenova\\
SCGP and Department of Mathematics, 3-115 \\
Stony Brook University \\
Stony Brook, NY 11794-3651, USA,} \\
\tt kamenova@math.stonybrook.edu
\\
 
	\noindent \sc Misha Verbitsky\\
		\sc Instituto Nacional de Matem\'atica Pura e
			Aplicada (IMPA) \\ Estrada Dona Castorina, 110\\
			Jardim Bot\^anico, CEP 22460-320
			Rio de Janeiro, RJ - Brasil \\
		also:\\
		Laboratory of Algebraic Geometry, \\
		Faculty of Mathematics, HSE University,\\
		6 Usacheva Str. Moscow, Russia\\
{	\tt verbit@impa.br }

\end{document}